\documentclass[12pt]{article}
\author{Gun Srijuntongsiri\thanks{4130 Upson Hall, Cornell University, Ithaca, NY 14853. Email: gunsri@cs.cornell.edu.},
Stephen A. Vavasis\thanks{MC 6054, University of Waterloo, 200
University Avenue W., Waterloo, ON N2L 3G1, Canada. Email:
vavasis@math.uwaterloo.ca.}}
\title{A Condition Number Analysis of a Line-Surface Intersection Algorithm\thanks{Supported in part by NSF DMS 0434338 and NSF CCF 0085969.}}
\date{\today}

\usepackage{amsmath, amsthm, amssymb,graphicx,subfig}

\begin{document}
\maketitle

\newcommand{\cond}[1]{\mathop{\rm{cond}}(#1)}
\newcommand{\norm}[1]{\left\|#1\right\|} %_\infty}
\newcommand{\norma}[1]{\left\|#1\right\|}   % arbitrary norm
\newcommand{\inv}[1]{#1^{-1}}
\newcommand{\nchoosek}[2]{\left(\begin{array}{c} #1\\ #2 \end{array} \right)}

\newtheorem{thm}{Theorem}[section]
\newtheorem{cor}[thm]{Corollary}
\newtheorem{lemma}[thm]{Lemma}
\newtheorem{prop}{Proposition}[thm]

\theoremstyle{remark}
\newtheorem{rem}[thm]{Remark}

\begin{abstract}
%We propose an algorithm based on Newton's method and subdivision for
%finding all zeros of a polynomial system in a bounded region of the
%plane.  This algorithm can be used to find the
%intersections between a line and a surface, which has
%applications in graphics and computer-aided geometric design.  Our
%analysis shows that the running time of the algorithm is bounded in
%terms of the condition number of the polynomial's zeros.  We argue that,
%in contrast, some other well known algorithms for this problem similarly
%based on Newton and subdivision may in certain cases require excessive
%computation even for well conditioned problems.
We propose an algorithm based on Newton's method and subdivision
for finding all zeros of a polynomial system in a bounded region
of the plane.  This algorithm can be used to find the
intersections between a line and a surface, which has applications
in graphics and computer-aided geometric design. The algorithm can
operate on polynomials represented in any basis that satisfies a
few conditions.  The power basis, the Bernstein basis, and the
first-kind Chebyshev basis are among those compatible with the
algorithm.  The main novelty of our algorithm is an analysis
showing that its running is bounded only in terms of the condition
number of the polynomial's zeros and a constant depending on the
polynomial basis.
\end{abstract}

\section{Introduction}

The problem of line-surface intersection has many applications in
areas such as geometric modeling, robotics, collision avoidance,
manufacturing simulation, scientific visualization, and computer
graphics. It is also a basis for considering intersections between
more complicated objects.  This article deals with intersections
of a line and a parametric surface.  The parametric method of
surface representation is a very convenient way of approximating
and designing curved surfaces, and computation using parametric
representation is often much more efficient than other types of
surface representations.

Typically, intersection problems reduce to solving systems of
nonlinear equations. Subdivision methods introduced by Whitted
\cite{whitted, rubi} were the first to be used for this problem.
In these methods, a simple shape such as rectangular box or sphere
is used to bound the surface and is tested whether the line
intersects the bounding volume.  If it does, the surface patch is
subdivided, and the bounding volumes are found for each subpatch.
The process repeats until no bounding volumes intersect the line
or the volumes are smaller than a specified size where the
intersection between such volumes and the line are taken as the
intersections between the surface and the line.  Subdivision
methods are robust and simple, but normally not efficient when
high accuracy of the solutions are required.  They also cannot
indicate if there are more than one zero inside the remaining
subdomains.

Regardless, a variation of subdivision methods known as B\'{e}zier clipping
by Nishita et al.\ should be noted for its efficient subdivision \cite{nishita}.
For a non-rational B\'{e}zier surface,
B\'{e}zier clipping uses the intersection between the convex hull of the
orthographic projection of the surface along the line
and a parameter axis to determine the regions which do not contain any intersections
before subdividing the remaining region. Sherbrooke and Patrikalakis
generalizes B\'{e}zier clipping to a zero-finding algorithm for an $n$-dimensional
nonlinear polynomial system called
\emph{Projected Polyhedron} (PP) algorithm \cite{sherbrooke}.

On the numerical side, Kajiya \cite{kajiya} proposes a method for intersecting
a line with a bicubic surface based on algebraic geometry without subdivisions.
His method is robust and simple. However, it is too costly to extend
to higher degree polynomials.
J\'{o}nsson and Vavasis \cite{jonsson} introduce
an algorithm for solving systems of two polynomials in two variables using
Macaulay resultant matrices.
They also analyze
the accuracy of the computed zeros in term of the conditioning of the problem.

Another approach is to combine a subdivision method with
a Newton-type method, using the latter to find the solutions of
the resulted system of equations after subpatches pass some criteria.
%Similar to pure subdivision methods,
%algorithms in this category subdivide the surface and use bounding volumes
%to eliminate the areas that cannot have intersections.  But when
%the remaining subpatch is small enough, an iterative method, usually Newton's method,
%is invoked to compute the exact intersection.
The advantages of Newton's method are
its quadratic convergence and simplicity in implementation, but it requires
a good initial approximation to converge and does not guarantee
that all zeros have been found.
To remedy these shortcomings, Toth \cite{toth}
uses the result from interval analysis to determine
the ``safe regions", where a Newton-like method starting from any point in them
are guaranteed to converge.  He tests each patch if it is a safe region
and if its axis-aligned bounding box intersects the line. If neither
is true, the patch is subdivided recursively.
Lischinski and Gonczarowski \cite{Lischinski}
propose an improvement to Toth's algorithm specific to
scene-rendering in computer graphics
by utilizing ray and surface coherences.

In contrast, other researchers develop
methods to estimate good initial points for Newton's method
rather than test for convergence of each choice of initial points.
These methods use tighter bounding volumes and subdivide the surface adaptively
until subpatches are \emph{flat} enough, that is, until they are close
enough to the bounding volumes.  Then the intersection between
the bounding volume and the line is chosen as the initial point
for Newton's method. Examples of methods in this category are
\cite{barth, fournier, qin, sweeney, yang}.
There is also the ray-tracing algorithm by Wang et al.\ that
combines Newton's method
and B\'{e}zier clipping together \cite{wang}.

Our algorithm is in the same category as Toth's in that it tests
for the convergence of initial points before performing Newton's
iteration to find solutions and it uses a bounding polygon of a
subpatch to exclude the one that cannot have a solution.  The
convergence test our algorithm uses is the Kantorovich test.
Because the Kantorovich test also tells us whether Newton's method
converges quadratically for the initial point in question in
addition to whether it converges at all, we can choose to hold off
Newton's method until quadratic convergence is assured.

The main feature of our algorithm is that there is an upper bound
on the number of subdivisions performed during the course of the
algorithm that depends only on the condition number of the problem
instance.  For example, having a solution located exactly on the
border of a subpatch does not adversely affect its efficiency.  To
the best of our knowledge, there is no previous algorithm in this
class whose running time has been bounded in terms of the
condition of the underlying problem instance, and we are not sure
whether such an analysis is possible for previous algorithms. Its
efficiency also depends on the choice of basis because the type of
bounding polygon depends on the basis.

The notion of bounding the running time of an iterative method
in terms of the condition number of the instance is an old one,
with the most notable example being the condition-number
bound of conjugate gradient (see Chapter 10 of \cite{gvl}).
This approach has also been used in interior-point methods
for linear programming \cite{freund} and Krylov-space eigenvalue
computation \cite{toh}.

\section{The theorem of Kantorovich}

Denote the closed ball centered at $x$ with radius $r>0$ by
\[
\bar{B}(x,r) = \{ y \in \mathbb{R}^n : \norm{y-x} \leq r \},
\]
and denote $B(x,r)$ as the interior of $\bar{B}(x,r)$.
Kantorovich's theorem in affine invariant
form is

\begin{thm}[Kantorovich, affine invariant form \cite{deuflhard, kantorovich}]
\label{standardkantorovich}
Let $f : D \subseteq \mathbb{R}^n \rightarrow \mathbb{R}^n$ be differentiable in
the open convex set $D$. Assume that for some point $x^0 \in D$, the Jacobian $f'(x^0)$
is invertible with
\[
\norm{f'(x^0)^{-1}f(x^0)} \leq \eta.
\]
Let there be a Lipschitz constant $\omega > 0$ for $f'(x^0)^{-1} f'$ such that
\[
\norm{f'(x^0)^{-1}(f'(x)-f'(y))} \leq \omega \cdot \norm{x-y} \textrm{ for all } x,y \in D.
\]
If $h = \eta\omega \leq 1/2$ and $\bar{B}(x^0,\rho_-) \subseteq D$, where
\[
\rho_- = \frac{1-\sqrt{1-2h}}{\omega},
\]
then $f$ has a zero $x^*$ in $\bar{B}(x^0,\rho_-)$. Moreover, this zero is the unique zero
of $f$ in $(\bar{B}(x^0,\rho_-) \cup B(x^0,\rho_+)) \cap D$ where
\[
\rho_+ = \frac{1+\sqrt{1-2h}}{\omega}
\]
and the Newton iterates $x^k$ with
\[
x^{k+1} = x^k - f'(x^k)^{-1}f(x^k)
\]
are well-defined, remain in $\bar{B}(x^0,\rho_-)$, and converge to $x^*$. In addition,
\begin{equation}
\label{rapidconv}
\norm{x^*-x^k} \leq \frac{\eta}{h}\left( \frac{(1-\sqrt{1-2h})^{2^k}}{2^k} \right), k = 0,1,2,\ldots
\end{equation}
\end{thm}
We call $x^0$ a \emph{fast starting point} if the sequence of Newton iterates starting from it
converges to a solution $x^*$ and (\ref{rapidconv}) is satisfied with $h \leq 1/4$,
which implies quadratic convergence of the iterates starting from $x^0$.
The Kantorovich's theorem also holds for complex functions \cite{farouki}.

\section{Formulation and representation of the line-surface
intersection problem} \label{section_formulation}

Let $\phi_0$, \ldots, $\phi_n$ denote a basis for the set of
univariate polynomials of degree at most $n$.  For example, the
power basis is defined by $\phi_i(t)=t^i$.  Other choices of basis
are discussed below.

Let $S$ be a two-dimensional surface embedded in $\mathbb{R}^3$
parametrized by
\[
\begin{array}{llll}
\bar{f}(u,v) &=& \sum_{i=0}^m \sum_{j=0}^n \bar{c}_{ij}
\phi_i(u)\phi_j(v), & 0 \leq u,v \leq 1,
\end{array}
\]
where $\bar{c}_{ij} \in \mathbb{R}^3$ $(i=0,1,\ldots,m;
j=0,1,\ldots,n)$ denote the coefficients. Define a line
\[
\begin{array}{llll}
L & =&  \{p+dt : & t \in \mathbb{R} \},
\end{array}
\]
where $p, d \in \mathbb{R}^3, d \neq 0$. The line-surface
intersection problem is to find all of the intersections between
$S$ and $L$, which are the solutions of the polynomial system
\begin{equation}
\label{inteq} \bar{f}(u,v)-(p+dt)=0.
\end{equation}
The system (\ref{inteq}) can be reduced to a system of two
equations and two unknowns. To show this, we first break
(\ref{inteq}) into its component parts
\begin{eqnarray}
\bar{f}_1(u,v)-p_1-td_1 & = &0, \nonumber \\
\bar{f}_2(u,v)-p_2-td_2 & = &0, \nonumber \\
\bar{f}_3(u,v)-p_3-td_3 & = &0. \nonumber
\end{eqnarray}
Here, the subscript $i$ denotes the $i$th coordinate of the point
in three-dimensional space. Assuming $|d_1| \geq \max
\left\{|d_2|,|d_3|\right\}$, we have the equivalent system
\begin{equation}
\left.
\begin{array}{lll}
d_1(\bar{f}_2(u,v)-p_2)-d_2(\bar{f}_1(u,v)-p_1) = 0, \\
d_1(\bar{f}_3(u,v)-p_3)-d_3(\bar{f}_1(u,v)-p_1) = 0,
\end{array} \right.
\label{red1}
\end{equation}
which can be rewritten with the same basis $\phi_i(u)\phi_j(v)$
(see item \ref{constant_poly_easy} on the list of basis properties below) as
\begin{equation}
\label{red2} f(u,v) \equiv \sum_{i=0}^m \sum_{j=0}^n c_{ij}
\phi_i(u)\phi_j(v) = 0.
\end{equation}
The system (\ref{red2}) is the one our algorithm operates on.

Since the parametric domain of the surface under consideration is
square, our algorithm uses the infinity norm for all of its norm
computation.  Therefore, for the rest of this article, the
notation $\norm{\cdot}$ is used to refer specifically to infinity
norm.

Our algorithm works with any polynomial basis $\phi_i(u)\phi_j(v)$
provided that the following properties hold:
\begin{enumerate}
\item There is a natural interval $[l,h]$ that is the domain for the polynomial.
In the case of Bernstein polynomials, this is $[0,1]$, and in the case of power and
Chebyshev polynomials, this is $[-1,1]$.

\item It is possible to compute a bounding polygon $P$ of $S =
\{f(u,v) : l \leq u,v \leq h \}$, where $f(u,v) = \sum_{i=0}^m
\sum_{j=0}^n c_{ij} \phi_i(u)\phi_j(v)$ and $c_{ij} \in
\mathbb{R}^d$ for any $d \geq 1$, that satisfies the following
properties: \label{bounding_prop}
\begin{enumerate}
\item Determining whether $0 \in P$ can be done efficiently
(ideally in $O(mn)$ operations).

\item Polygon $P$ is affinely and translationally invariant. In other words,
the bounding polygon of $\{Af(u,v)+b: l \leq u,v \leq h\}$
 is $\{Ax+b : x \in P \}$ for any nonsingular matrix
$A \in \mathbb{R}^{d \times d}$ and any vector $b \in \mathbb{R}^d$.

\item For any $y \in P$,
\begin{equation}
\norm{y} \leq \theta\max_{l \leq u,v \leq h} \norm{ f(u,v) },
\label{thetadef}
\end{equation}
where $\theta$ is a function of $m$ and $n$.
\item If $d=1$, then the endpoints of $P$ can be computed efficiently (ideally in $O(mn)$ time).
\label{prop_endpoint}
\end{enumerate}

\item It is possible to reparametrize with $[l,h]^2$ the surface
$S_1 = \{ f(x) : x \in \bar{B}(x^0,r) \}$, where $x^0 \in
\mathbb{R}^2$ and $r \in \mathbb{R} > 0$. In other words, it is
possible (and efficient) to compute the polynomial $\hat f$
represented in the same basis such that $S_1=\{\hat f(\hat x):
\hat x \in [l,h]^2\}$.

\item Constant polynomials are easy to represent. \label{constant_poly_easy}

\item Derivatives of polynomials are easy to determine in the same basis. (preferably in $O(mn)$ operations).
\label{prop_deriv}
\end{enumerate}
Recall that $P$
is a bounding polygon of $S$ if and only if $x \in S$ implies $x
\in P$.

As shown later in Section \ref{section_analysis}, the efficiency
of our algorithm depends on $\theta$. Hence, the choice
of the basis affects the algorithm's performance as each basis
allows different ways to compute bounding polygons.
%The readers
%are referred to \cite{srijuntongsiri2} for detail analysis of
%$\theta$ and $\tau$ for power, Bernstein, and Chebyshev bases.

\section{The Kantorovich-Test Subdivision algorithm}

Before we detail our algorithm, we define notation and crucial
quantities that are used by the algorithm and its analysis.
Denote $x = (u,v)$ as a point in two-dimensional parametric space
and $f(x) = f(u,v)$ as the value of $f$ at $x$.

For a given zero $x^*$ of polynomial $f$,
let $\omega_*(x^*)$ and $\rho_*(x^*)$
be quantities satisfying the conditions that, first, $\omega_*(x^*)$ is the smallest Lipschitz constant
for $f'(x^*)^{-1}f'$, i.e.,
\begin{equation}
\label{lip1}
\norm{f'(x^*)^{-1}\left( f'(x)-f'(y) \right)} \leq \omega_*(x^*) \cdot \norm{x-y} \textrm{ for all } x,y \in \bar{B}(x^*,\rho_*(x^*))
\end{equation}
and, second,
\begin{equation}
\label{rhostar}
\rho_*(x^*) = \frac{2}{\omega_*(x^*)}.
\end{equation}
Since $\omega_*(x^*)$ is nondecreasing as $\rho_*(x^*)$ increases in (\ref{lip1}) but
$\rho_*(x^*)$ is decreasing as $\omega_*(x^*)$ increases in (\ref{rhostar}), there exists a unique pair
$(\omega_*(x^*),\rho_*(x^*))$ satisfying the above conditions, and this pair can be obtained by binary search.
When more than one function is being discussed, we use $\omega_*^f(x^*)$ to denote $\omega_*(x^*)$ of
the function $f$. Approximation to these two quantities, $\omega_*(x^*)$ and $\rho_*(x^*)$, are computed and made use of by
the algorithm.

For clarity, we simply abbreviate
$\omega_*(x^*)$ and $\rho_*(x^*)$ as $\omega_*$ and $\rho_*$, respectively,
throughout the rest of this article when it is clear from the context to which
$x^*$ the quantities belong.

Straightforward application of the affine invariant form of Kantorovich's theorem with $x^0=x^*$ and $D = \bar{B}(x^*,\rho_*(x^*))$
yields the result that $x^*$ is the unique zero of $f$ in $\bar{B}(x^*,\rho_*(x^*))$.  In fact,
the above definitions of $\omega_*(x^*)$ and
$\rho_*(x^*)$ are chosen such that the ball that is guaranteed by the
Kantorovich's theorem to contain no other zeros than $x^*$ is the largest possible.

Define
\[
\gamma(\theta) = 1 / \left(4\sqrt{\theta(4\theta+1)}-8\theta \right),
\]
where $\theta$ is as in (\ref{thetadef}).  Observe that $1 \leq \gamma(\theta) \leq \left(\sqrt{5}+2\right)/4 \approx 1.0590$ since $\gamma(\theta)$ is a decreasing function for positive $\theta$ and $\theta \geq 1$ by
the definition of a bounding polygon.
Another quantity of interest is $\omega_{D'}$, which is defined as
the smallest nonnegative constant $\omega$ satisfying
\begin{equation}
\begin{array}{llll}
\norm{f'(x^*)^{-1}\left(f'(y)-f'(z)\right)} & \leq & \omega \cdot \norm{y-z}, & y,z \in D', x^* \in [0,1]^2
 \\
& & & \textrm{satisfying } f(x^*) = 0,  \label{om}\end{array}
\end{equation}
where
\begin{equation}
\label{dprimedef} D' =
\left[-\gamma(\theta),1+\gamma(\theta)\right]^2.
\end{equation}
The motivation of this definition of $D'$ is that it contains all
domains whose Lipschitz constants may be needed during the course
of the algorithm. Denote $\omega_f$ as the maximum of
$\omega_{D'}$ and all $\omega_*(x^*)$
\[
\omega_f = \max\{\omega_{D'}, \max_{x^* \in \mathbb{C}^2 : f(x^*)=0} \omega_*(x^*)\}.
\]
Finally, define the condition number of $f$ to be
\begin{equation}
\label{conddef}
\cond{f} = \max \{ \omega_f, \max_{x^* \in \mathbb{C}^2 : f(x^*)=0,y \in [0,1]^2} \norm{\inv{f'(x^*)}f'(y)}\}.
\end{equation}
Note that in (\ref{om}), $x^*$ is restricted to zeros in $[0,1]^2$ whereas in (\ref{conddef}),
$x^*$ ranges over all complex zeros of $f$.
We defer the discussion of why (\ref{conddef}) is a reasonable condition number until
after the description of our algorithm.

We define the \emph{Kantorovich test} on a region $X =
\bar{B}(x^0,r)$ as the application of Kantorovich's Theorem on
the point $x^0$ using $\bar{B}(x^0,2\gamma(\theta)r)$ as the
domain $D$ in the statement of the theorem and
$\norm{f'(x^0)^{-1}f(x^0)}$ as $\eta$. For $\omega$, we instead
use $\hat{\omega} \geq \omega$, where $\hat{\omega}$ is defined by
(\ref{computedlip}) below. The region $X$ passes the Kantorovich
test if $\eta \hat{\omega} \leq 1/4$ and $\bar{B}(x^0,\rho_-)
\subseteq D'$, which implies that $x^0$ is a fast starting point.

The other test our algorithm uses is the exclusion test. For a
given region $X$, let $\hat{f}_X$ be the polynomial in the basis
$\phi_i(u)\phi_j(v)$ that reparametrizes with $[l,h]^2$ the
surface defined by $f$ over $X$.  The region $X$ passes the
\emph{exclusion test} if the bounding polygon of $\{
\hat{f}_X(u,v) : l \leq u,v \leq h \}$ excludes the origin.  Note
that the bounding polygon used in this test must satisfy
item \ref{bounding_prop} of the basis properties listed in Section \ref{section_formulation}.
%The
%choice of which type of bounding polygons to use depends on the
%choice of the polynomial basis.  A suitable bounding polygon
%should be easy to compute and determine if it intersects with
%another polygon of the same type. For most bases, an axis-aligned
%bounding boxes can be used as the bounding polygon.  For Bernstein
%basis, a convex hull is a better choice as it can be computed
%directly from the coefficients and it bounds the surface tighter
%than the axis-aligned bounding box.  Note that computing a convex
%hull of a surface represented by other basis is generally too
%expensive to compensate for the tighter bound it offers over other
%kinds of bounding polygons.

We now proceed to describe our algorithm, the
\emph{Kantorovich-Test Subdivision algorithm} or KTS in short.

\begin{flushleft}
\textbf{Algorithm KTS}:
\end{flushleft}
\begin{itemize}
\item Let $Q$ be a queue with $[0,1]^2$ as its only entry. Set $S = \emptyset$.
\item Repeat until $Q = \emptyset$
\begin{enumerate}
\item Let $X$ be the patch at the front of $Q$.  Remove $X$ from $Q$.
\item If $X \not\subseteq X_S$ for all $X_S \in S$,
\begin{itemize}
\item Perform the exclusion test on $X=\bar{B}(x^0,r)$
\item If $X$ fails the exclusion test,
\begin{enumerate}
\item Perform the Kantorovich test on $X$
\item If $X$ passes the Kantorovich test,
\begin{enumerate}
\item Perform Newton's method starting from $x^0$ to find a zero $x^*$.
\item If $x^* \not\in X_S$ for any $X_S \in S$ (i.e., $x^*$ has not been found previously),
\begin{itemize}
\item Compute $\rho_*(x^*)$ and its associated $\omega_*(x^*)$ by binary search.
\item Set $S = S \cup \left\{\bar{B}(x^*,\rho_*(x^*))\right\}$.
\end{itemize}
\end{enumerate}
\item Subdivide $X$ along both $u$ and $v$-axes into four equal subregions. Add these subregions
to the end of $Q$.
\end{enumerate}
\end{itemize}
\end{enumerate}
\end{itemize}
A few remarks are needed regarding the description of the KTS algorithm.
\begin{itemize}
\item The subdivision in step 2.c is performed regardless of the result of the Kantorovich test.
In general, passing the Kantorovich test does not imply that there is only one zero in $X$.
\item The check that the zero found
by Newton's method is not a duplicate (step 2.b.ii) is necessary
since the Kantorovich test may detect a zero outside $X$.
\item If the Kantorovich test is not applicable
for a certain patch due to the Jacobian of the midpoint being singular, the patch is
treated as if it fails the Kantorovich test.
\end{itemize}

One property of KTS is that it is affine invariant.  In other words, left-multiplying $f$
with a 2-by-2 matrix $A$ prior to executing KTS does not change its behavior.
This is the main reason we define the condition number to be affine invariant.
Define $g \equiv Af$.
To see that our condition number is affine invariant, note that
$\inv{g'(x)}g'(y) = \inv{[Af'(x)]}Af'(y) = \inv{f'(x)}\inv{A}Af'(y) = \inv{f'(x)}f'(y)$
for any $x, y \in \mathbb{R}^n$. Therefore, $\cond{g} = \cond{f}$.
In contrast, simpler condition numbers such as Lipschitz constants for $f'$
are not affine invariant and hence are not chosen for our analysis.

Since Toth's algorithm is the most similar one to KTS, it is worthwhile to discuss
the main differences between the two and the implications these differences make.  First,
Toth's uses the Krawczyk-Moore test and another unnamed test, both based on interval
analysis, as the convergence test.  These two tests guarantee linear convergence for the \emph{simple} Newton iteration---a variation of Newton's method where the Jacobian of the initial
point is used in place of the Jacobian of the current point in every iteration.
With our Kantorovich test, KTS starts Newton's method only when quadratic
convergence is assured.

Another main difference is in the choice of domains for the convergence test.
Toth's uses the subpatch $X$ itself as the domain for the test.  This choice
may exhibit undesirable behavior when a zero lies on the border of a subpatch,
which is not necessarily on or near the border of the original domain $[0,1]^2$.
For example, consider the function $f(u,v)=(u^2-.25,v-.8)^T$ whose zeros are $(.5,.8)$
and $(-.5,.8)$.  The patch $\{ (u,v) : .5 \leq u \leq .5+\epsilon, a \leq v \leq b \}$
does not pass either of Toth's convergence tests for any $\epsilon > 0$ and any $a \leq .8 \leq b$
although the patch $\{ (u,v) : .45 \leq u \leq .8, 0 \leq v \leq 1 \}$, a large patch whose
borders do not coincide with any zeros, does pass the Krawczyk-Moore test.  This results
in excessive subdivisions by Toth's algorithm.  KTS uses $\bar{B}(x^0,2\gamma(\theta)r)$ as
the domain for $X$ to avoid this problem.  Theorem \ref{thm1} below shows that the Kantorovich test
does not have trouble detecting the zeros locating on the border of the subpatch.

\section{Implementation details when using power, Bernstein, or Chebyshev bases}

This section covers the implementation details of KTS when the
polynomial system is in the power, Bernstein, or Chebyshev
bases.  The power basis for polynomials of degree $n$ is $\phi_k(t) =
t^k$ $(0 \leq k \leq n)$. The Bernstein basis is $\phi_k(t) =
Z_{k,n}(t) = \nchoosek{n}{k}(1-t)^{n-k}t^k$ $(0 \leq k \leq n)$.
The Chebyshev basis is $\phi_k(t) = T_k(t)$  $(0 \leq k \leq n)$,
where $T_k(t)$ is the Chebyshev polynomial of the first kind
generated by the recurrence relation
\begin{eqnarray}
T_0(t) &= &1, \nonumber \\
T_1(t) &= & t, \nonumber \\
T_{k+1}(t) & = & 2tT_k(t)-T_{k-1}(t) \textrm{ for } k \geq 1.
\label{cheb_recurrence}
\end{eqnarray}

\subsection{Bounding polygons}
\label{section_bounding_polygon} We begin with the choices of $l$
and $h$ and the definitions of bounding polygons of the surface
$S = \{ f(u,v) : l \leq u,v \leq h\}$, where $f(u,v)$ is represented
by one of the three bases, that satisfy the required properties
detailed in Section \ref{section_formulation}. For Bernstein
basis, the convex hull of the coefficients (control points), call
it $P_1$, satisfies the requirements for $l = 0$ and $h=1$.  The
convex hull $P_1$ can be described as
\begin{equation}
\label{p1def} P_1 = \left\{ \sum_{i,j} c_{ij} s_{ij} : \sum_{i,j}
s_{ij} = 1, 0 \leq s_{ij} \leq 1 \right\}.
\end{equation}
For power and Chebyshev bases, the bounding polygon
\begin{equation}
\label{p2def} P_2 = \left\{ c_{00} + \sum_{i+j > 0} c_{ij} s_{ij} : -1
\leq s_{ij} \leq 1 \right\}
\end{equation}
satisfies the requirements for $l = -1$ and $h = 1$.  Note that
$P_2$ is a bounding polygon of $S$ in the Chebyshev case since $|T_k(t)| \leq 1$
for any $k \geq 0$ and any $t \in [-1,1]$.
Determining whether $0 \in P_2$ is done by solving a small linear programming problem.
The value of $\theta$ for each case is summarized in Table
\ref{table_theta}.  The reader is referred to \cite{srijuntongsiri2} for
the derivation of $\theta$ for each of the three bases
as well as the proofs that $P_1$ and $P_2$ satisfy all of the basis properties listed in Section \ref{section_formulation}.

\begin{table}
\begin{center}
\begin{tabular}{|c|c|}
\hline
Basis & $\theta$ \\
\hline \hline
%Bernstein & $0$ & $1$ & convex hull of coefficients & $(m+1)(n+1)m^{m} n^n$ \\
Bernstein &  %convex hull of coefficients &
$\left(\sum_{i=0}^m \prod_{i' \neq i} \frac{\max\{|m-i'|,|i'|
\}}{|i-i'|} \right)
\left(\sum_{j=0}^n \prod_{j' \neq j} \frac{\max\{|n-j'|,|j'| \}}{|j-j'|} \right) = O(m^{m+1}n^{n+1})$ \\
Chebyshev & %$\{ c_{00} + \sum_{i+j > 0} c_{ij} s_{ij} : -1 \leq s_{ij} \leq 1 \}$
$2(m+1)(n+1)$ \\
Power & %$\{ c_{00} + \sum_{i+j > 0} c_{ij} s_{ij} : -1 \leq s_{ij} \leq 1 \}$
$(m+1)(n+1)(3^{m+1}-1)(3^{n+1}-1)/2$ \\
\hline
\end{tabular}
\end{center}
\caption{The value of $\theta$'s of the power,
the Bernstein, and the Chebyshev bases and their corresponding bounding polygons.} \label{table_theta}
\end{table}

\subsection{Computation of Lipschitz constant}
\label{section_imp} Another step of KTS that needs further
elaboration is the computation of Lipschitz constant in the
Kantorovich test.  The Lipschitz constant for $f'(x^0)^{-1}f'
\equiv g$ is obtained from an upper bound of the derivative of $g$
\[
g'(x) =
\left(\frac{\partial^2\left(f'(x^0)^{-1}f\right)_i(x)}{\partial
x_j
\partial x_k} \right)
\]
for all $x \in X$.  Let $\hat{g} \equiv \hat{g}_X$ be the
polynomial in the same basis as $f$ that reparametrizes with
$[l,h]^2$ the surface defined by $g$ over $X$. We have that
\begin{eqnarray}
\max_{x \in X} \norm{g'(x)} & = & \max_{x \in [l,h]^2} \norm{\hat{g}'(x)} \nonumber \\
& = & \max_{x \in [l,h]^2} \max_{\norm{y}=1} \norm{\hat{g}'(x)y} \nonumber \\
                            & \leq & \max_{x \in [l,h]^2} \max_i \sum_{j=1}^2\sum_{k=1}^2 |\hat{g}'_{ijk}(x)| \nonumber \\
                            & \leq & 4 \max_{i,j,k} \max_{x \in [l,h]^2} | \hat{g}'_{ijk}(x)
                            |.\nonumber
%                           & \leq & 4 \max_{i,j,k} \max_{x \in X} U_{ijk}, \nonumber
\end{eqnarray}
Note that each entry of $\hat{g}'$ can be written as a polynomial
in the same basis as $f$ (refer to property \ref{prop_deriv} of the basis).  For this reason, an upper bound of
$\max_{x \in [l,h]^2} | \hat{g}'_{ijk}(x)|$ can be computed as
follows: Let $P_{ijk}$ be the bounding polygon (bounding interval
in this case) of $\{\hat g'_{ijk}(x): x\in[l,h]^2\}$
computed in the same way as described in Section
\ref{section_bounding_polygon}. The maximum absolute value of the
endpoints of $P_{ijk}$ (refer to property \ref{prop_endpoint} of the basis) is an upper bound of $\max_{x \in [l,h]^2}
| \hat{g}'_{ijk}(x)|$. Let $\hat{\omega}$ denote the Lipschitz
constant computed in this manner, that is,
\begin{equation}
\label{computedlip} \hat{\omega} \equiv 4 \max_{i,j,k} \max_{x \in [l,h]^2} | \hat{g}'_{ijk}(x)|,
\end{equation}
where $\max_{x \in [l,h]^2} | \hat{g}'_{ijk}(x)|$ is computed from the
endpoints of its bounding interval.

\section{Significance of our condition number}

We now discuss the significance of (\ref{conddef}) to the
conditioning of the problem. In particular, we attempt to justify
that the efficiency of any algorithm in the same class as KTS is
dependent on (\ref{conddef}). This class of algorithms being
considered includes any algorithm that (i) isolates unique zeros
with subdivision before finding them and (ii) will not discard a
patch until the convex hull of its function values (which is
clearly a subset of any possible bounding convex polygon) excludes
the origin.

\subsection{Condition number and the Kantorovich test}

This section discusses the relationship between $\omega_f$ and the
Kantorovich test.  We show that, for any given zero $x^*$ of an
arbitrary $f$, there is a function $\bar{f}$ such that $x^*$ is
also a zero of $\bar{f}$, $f'(x^*) = \bar{f}'(x^*)$,
$\omega_*^f(x^*) = \omega_*^{\bar{f}}(x^*)$, and $\bar{f}$ has
another zero $y^*$ with $\norm{y^*-x^*} = \rho_*(x^*)$.
For example, consider a zero $x^* = (.5,.5)$ of
the function $f = (u^3-2.2 u^2+1.55 u-.35,v^2-.7 v+.1)^T$,
of which $\rho_*(x^*) = .1$.
A corresponding $\bar{f}$ with the above properties is
$\bar{f} = (u^2-.9 u+.2, .3 v -.15)^T$,
which has zeros at $(.5,.5)$ and $(.4,.5)$.
Since the
Kantorovich test uses only the function value, its first
derivative, and the Lipschitz constant, all of which are the same
for $f$ and $\bar{f}$ at $x^*$, the functions $f$ and $\bar{f}$
are identical from the perspective of the Kantorovich test applied
to $x^*$. Therefore, $\rho_*(x^*)$ is a reasonable number that
quantifies the distance between $x^*$ and its nearest other zero
barring the usage of additional information. Consequently,
$\omega_f$, which is greater than or equal to $\omega_*(x^*) =
2/\rho_*(x^*)$ for all zeros $x^*$ of $f$, describes the distance
between the closest pairs of zeros of $f$. Therefore, the
efficiency of any algorithm that isolates unique zeros is
dependent on $\omega_f$.
%\begin{figure}
%% \centering
% \subfloat[The contour lines of $f$ and the zeros of $f$. The zero $x^*$ has $\omega_*^f(x^*)=20$ and $\rho_*(x^*)=.1$.]{ \includegraphics[width=0.9\textwidth]{kanto_fig_f}}\qquad
% \subfloat[The contour lines of $\bar{f}$ and the two zeros $x^*=(.5,.5)$ and $y^*=(.4,.5)$ of $\bar{f}$. The distance between the two zeros $\norm{y^*-x^*} = .1 = \rho_*(x^*)$.]{\includegraphics[width=0.9\textwidth]{kanto_fig_fbar}\label{f1b}}
% \caption{The function $f = (u^3-2.2 u^2+1.55 u-.35,v^2-.7 v+.1)^T$ and its corresponding
% $\bar{f} = (u^2-.9 u+.2, .3 v -.15)^T$.  Choose $x^* = (.5,.5)$. The two functions $f$ and $\bar{f}$ satisfies $f(x^*)=\bar{f}(x^*)=0$, $f'(x^*)=\bar{f}'(x^*)$, $\omega_*^f(x^*)=\omega_*^{\bar{f}}(x^*)$, and $\bar{f}$ has
% another zero $y^*$ with $\norm{y^*-x^*} = \rho_*(x^*)$.}
% \label{fig_kanto}
%\end{figure}
%

The function $\bar{f}$ with the above properties can be
constructed as follows: Let $x^* = (u^*,v^*),$ $f'(x^*) = \left( \begin{array}{cc}
\alpha_1 & \alpha_2 \\ \alpha_3 & \alpha_4 \end{array} \right)$,
%[\alpha_1, \alpha_2 ; \alpha_3, \alpha_4]$,
and $\omega_*^f(x^*) =
\omega$. If $|\alpha_4| \geq |\alpha_3|$,
\begin{equation}
\label{ggg1}
\bar{f}(u,v) =          \left( \begin{array}{c}
                        \frac{\omega(\alpha_1\alpha_4-\alpha_2\alpha_3)}{2\alpha_4}(u-u^*)^2+\alpha_1(u-u^*)+\alpha_2 (v-v^*)  \\
                        \alpha_3 (u-u^*)  +\alpha_4 (v-v^*)\end{array} \right).
\end{equation}
Otherwise,
\[
\bar{f}(u,v) =          \left( \begin{array}{c}
                        \alpha_1(u-u^*)+\frac{\omega(\alpha_1\alpha_4-\alpha_2\alpha_3)}{2\alpha_3}(v-v^*)^2+\alpha_2(v-v^*) \\
                        \alpha_3 (u-u^*)  +\alpha_4 (v-v^*)
                        \end{array} \right).
\]

It is straightforward to verify that $\bar{f}(x^*) = 0$, $\bar{f}'(x^*) = f'(x^*)$, and
$\omega_*^{\bar{f}}(x^*) = \omega$.
We now show that $\norm{y^*-x^*} = \rho_*(x^*)$ for the
case where $|\alpha_4| \geq |\alpha_3|$.
The other case can be verified in the same manner.  Let $y^* = (u^* + \Delta u, v^*+\Delta v)$.
Substituting $y^*$ into (\ref{ggg1}) and setting it to zero yields
\begin{equation}
\label{gsimple}
g(u^* + \Delta u, v^*+\Delta v) = \left( \begin{array}{c}
                                         \frac{\omega(\alpha_1\alpha_4-\alpha_2\alpha_3)}{2\alpha_4}(\Delta u)^2 + \alpha_1 \Delta u + \alpha_2 \Delta v \\
                                         \alpha_3 \Delta u + \alpha_4 \Delta v \end{array} \right)
                                = \left( \begin{array}{c} 0\\ 0 \end{array} \right).
\end{equation}
Solving (\ref{gsimple}) yields
\begin{eqnarray}
\Delta u & = & -\frac{2}{\omega}, \nonumber \\
\Delta v & = & \frac{\alpha_3}{\alpha_4}\cdot\frac{2}{\omega}. \nonumber
\end{eqnarray}
Since $|\alpha_4| \geq |\alpha_3|$ and $\rho_*(x^*) = 2/\omega$, $\norm{y^*-x^*} = \rho_*(x^*)$.

\subsection{Condition number and the exclusion test}

The other term in our condition number, $\max_{x^* \in
\mathbb{C}^2 : f(x^*) = 0, y \in [0,1]^2}
\norm{\inv{f'(x^*)}f'(y)}$, relates to the convex bounding polygon
test---the test to determine whether the convex bounding polygon
of a subpatch contains the origin. We show that there exists a
function $f$ such that a patch $B(x^0,r)$ where $x^0$ is
relatively close to a zero, fails the convex bounding polygon test
if $r \geq O(1/\cond{f})$. Denote $x^0 = (u^0,v^0)$.  Define the
complex function $g(z) =
\left(z-(u^0-\epsilon-i\epsilon)\right)\cdot\left(z-(u^0+\epsilon-i\epsilon)\right)$,
where $\epsilon \in \mathbb{R}$ and $0 < \epsilon < 1$. Consider
the following function
\begin{eqnarray}
f(u,v) & =& \left( \begin{array}{l}
%                   u^2-v^2-2u^0u-2\epsilon v - 2\epsilon^2 + (u^0)^2\\
%                   2uv-2u^0v+2\epsilon u - 2 \epsilon u^0
                    {\rm Re}\left( g(u+iv) \right)\\
                    {\rm Im}\left( g(u+iv) \right)
                \end{array} \right)  \nonumber \\
                & = & \left( \begin{array}{l}
                   u^2-v^2-2u^0u-2\epsilon v - 2\epsilon^2 + (u^0)^2\\
                   2uv-2u^0v+2\epsilon u - 2 \epsilon u^0
                \end{array} \right), \label{fexc}
\end{eqnarray}
%where $\Re \epsilon > 0$.
where ${\rm Re}(z)$ and ${\rm Im}(z)$ denote the real and imaginary parts of the complex number $z$, respectively.
The four complex zeros of $f$ are $(u^0-\epsilon,-\epsilon)$, $(u^0+\epsilon,-\epsilon)$,
$(u^0,-\epsilon-i\epsilon )$, and $(u^0,-\epsilon+i\epsilon)$. Therefore,
\[
\max_{x^* \in \mathbb{C}^2 : f(x^*) = 0, y \in [0,1]^2} \norm{\inv{f'(x^*)}f'(y)} = O(1/\epsilon).
\]
Moreover,
\[
\omega_f = O(1/\epsilon).
\]

\begin{figure}
% \centering
 \subfloat[The circular arc $A \subseteq \bar{B}(x^0,r)$ .]{ \includegraphics[width=0.9\textwidth]{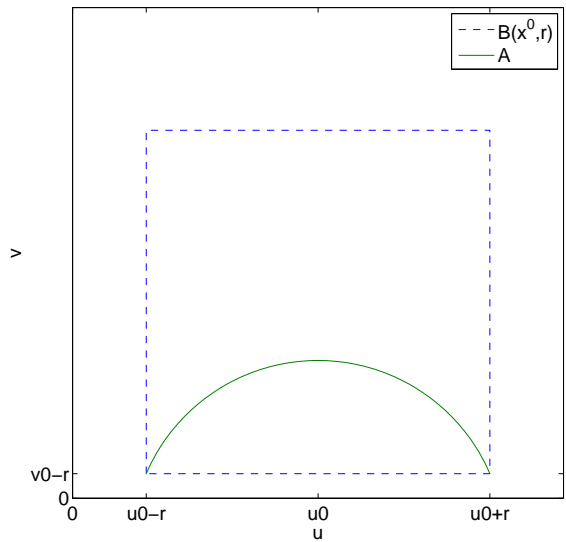}}\qquad
 \subfloat[The range $f(A)$ and its convex bounding polygon]{\includegraphics[width=0.9\textwidth]{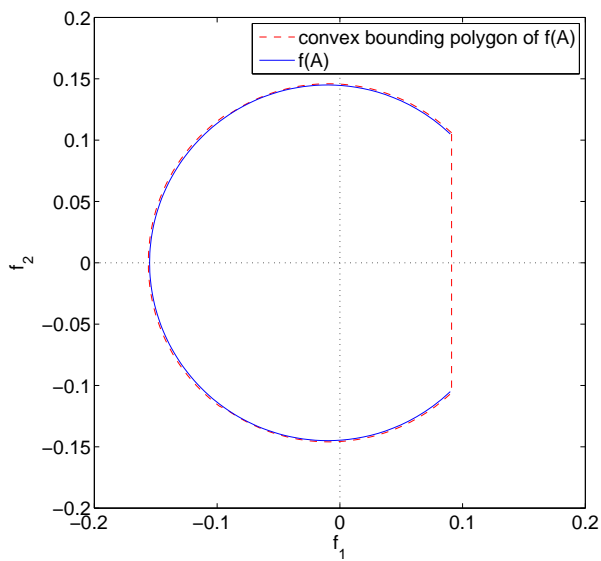}\label{f2b}}
 \caption{The circular arc $A$ centered at $(u^0,-\epsilon)$
 that goes from $(u^0+r,v^0-r)$ to $(u^0-r,v^0-r)$ and its range $f(A)$ where $f$ is
 as in (\ref{fexc}). Figure \ref{f2b} shows that the bounding convex
 polygon of $f(A)$ contains the origin, and therefore $\bar{B}(x^0,r)$
 fails the convex bounding polygon test.}
 \label{fig_exclusion}
\end{figure}

We now show for the case that $v^0 = O(\epsilon)$ that $B(x^0,r)$
fails the convex bounding polygon test if $r \geq O(\epsilon)$.
Let $A$ be the circular arc centered at $(u^0,-\epsilon)$ that
goes from $(u^0+r,v^0-r)$ to $(u^0-r,v^0-r)$ counterclockwise.
Observe that $f$ maps $A$ to the circular arc centered at
$(-\epsilon^2,0)$ that goes from $(2(v^0+\epsilon)r-2\epsilon
v^0-(v^0)^2-2\epsilon^2,2r(v^0+\epsilon-r))$ to
$(2(v^0+\epsilon)r-2\epsilon
v^0-(v^0)^2-2\epsilon^2,-2r(v^0+\epsilon-r))$ counterclockwise
(see Figure \ref{fig_exclusion}).
Notice that $2r(v^0+\epsilon-r) \geq 0$ because $B(x^0,r)
\subseteq [0,1]^2$. Therefore, the convex bounding polygon of
$f(A)$ contains the origin if $r > ((v^0)^2+2\epsilon v^0 +
2\epsilon^2) / (2(v^0+\epsilon)) = O(\epsilon)$ (recall the
assumption that $v^0 = O(\epsilon)$). Since $A \subset B(x^0,r)$,
the convex bounding polygon of $f(B(x^0,r))$ also contains the
origin and the convex bounding polygon test fails.
\section{Time complexity analysis}
\label{section_analysis}

In this section, we prove a number of theorems relating to the
behavior of the KTS algorithm.  We analyze the efficiency of KTS
by showing that a patch either is a subset of a safe region,
passes the Kantorovich test, or passes the exclusion test when it
is smaller than a certain size that depends on the condition
number of the function. Hence, we have the upper bound of the
total number of patches examined by KTS in order to solve the
intersection problem.

Recall that the Lipschitz constant $\hat{\omega}$ given by
(\ref{computedlip}) is not the smallest Lipschitz constant of
$\inv{f'(x^0)}f$ over $D'$, where $D'$ is given by
(\ref{dprimedef}). However, we can show that $\hat{\omega} \leq
4\theta\omega$, where $\omega$ denotes the smallest Lipschitz
constant of $\inv{f'(x^0)}f$ over $D'$. Since $\hat{\omega}$
is computed from the
endpoints of the bounding intervals of $\max_{x \in [l,h]^2} | \hat{g}'_{ijk}(x)|$,
by (\ref{thetadef}),
\begin{eqnarray}
\hat{\omega} & \leq & 4\theta \max_{i,j,k}\max_{x \in [l,h]^2} \left|\hat{g}'_{ijk}(x) \right| \nonumber \\
& = & 4\theta \max_{i,j,k}\max_{x \in X} \left|g'_{ijk}(x) \right| \nonumber \\
& \leq & 4\theta \max_{x \in X} \norm{g'(x)} = 4\theta\omega.
\label{hodef}
\end{eqnarray}
With this bound on $\hat{\omega}$, we can now analyze the behavior of the Kantorovich test.

\begin{thm}
\label{thm1} Let $x^0$ be a point in $[0,1]^2$ such that $f'(x^0)$
is invertible. Let $x^*$ be a zero of $f$ that is contained in
$\bar{B}(x^0,r)$, where $r$ is the radius of the patch under
consideration. The patch $X = \bar{B}(x^0,r) \subseteq [0,1]^2$
passes the Kantorovich test if
\begin{equation}
\label{thm1a1}
r \leq \frac{\gamma(\theta)-1}{\gamma(\theta)\omega_{D'}}.
\end{equation}
\end{thm}
\begin{proof}
The first step is to show that $\eta\hat{\omega} \leq 1/4$, where
$\hat{\omega}$ is as in (\ref{computedlip}). Since $r \leq 1/2$,
$\bar{B}(x^0,2\gamma(\theta)r) \subseteq D'$. Observe that for any
$x, y \in D'$,
\begin{eqnarray}
\norm{\inv{f'(x^0)}(f'(x)-f'(y))} & = & \|\left(\inv{f'(x^*)} + (\inv{f'(x^0)}-\inv{f'(x^*)})\right)
\nonumber \\
& & \left(f'(x)-f'(y)\right)\| \nonumber \\
& \leq & \norm{\inv{f'(x^*)}\left(f'(x)-f'(y)\right)}+ \|\inv{f'(x^*)} \nonumber \\
& & (f'(x^*)-f'(x^0))\inv{f'(x^0)}(f'(x)-f'(y))\| \nonumber \\
& \leq & \omega_{D'} \norm{x-y} + \norm{\inv{f'(x^*)}(f'(x^*)-f'(x^0))}\cdot \nonumber \\
& & \norm{\inv{f'(x^0)}(f'(x)-f'(y))} \nonumber \\
& \leq & \omega_{D'} \norm{x-y} + \nonumber \\
& &  \omega_{D'}\norm{x^*-x^0} \cdot \norm{\inv{f'(x^0)}(f'(x)-f'(y))} \nonumber\\
& \leq & \omega_{D'} \norm{x-y} + \nonumber \\
& & \omega_{D'} r \cdot \norm{\inv{f'(x^0)}(f'(x)-f'(y))}.
\label{thm1l1}
\end{eqnarray}
Since (\ref{thm1a1}) implies
\begin{equation}
\label{onewr} 1-\omega_{D'}r \geq 1/\gamma(\theta) > 0,
\end{equation}
the inequality (\ref{thm1l1}) becomes
\[
\norm{\inv{f'(x^0)}(f'(x)-f'(y))} \leq \left(\frac{\omega_{D'}}{1-\omega_{D'} r}\right) \norm{x-y}.
\]
Hence
\begin{equation}
\omega \leq \frac{\omega_{D'}}{1-\omega_{D'} r}, \label{w}
\end{equation}
where $\omega$ is the smallest Lipschitz constant of
$\inv{f'(x^0)}f'$ over $D'$.

Recall that $f(x^*) = 0$ and $X \subseteq D'$. Observe that
\begin{eqnarray}
\eta & \equiv & \norm{f'(x^0)^{-1}f(x^0)}  \nonumber \\
& = & \norm{f'(x^0)^{-1}(f(x^0)-f(x^*))} \nonumber \\
& \leq & \left(\max_{x \in X} \norm{\inv{f'(x^0)}f'(x)}\right)\cdot \norm{x^0-x^*} \nonumber \\
& \leq & \left(\max_{x \in X} \norm{\inv{f'(x^0)}(f'(x)-f'(x^0))+\inv{f'(x^0)}f'(x^0)}\right)\cdot r \nonumber \\
& \leq & \left(\max_{x \in X} \norm{\inv{f'(x^0)}(f'(x)-f'(x^0))}+1\right)\cdot r \nonumber  \\
& \leq & (\omega r + 1 ) r. \label{et}
\end{eqnarray}

Using (\ref{hodef}), (\ref{thm1a1}), (\ref{onewr}), (\ref{w}), and
(\ref{et}) together give
\[
\eta\hat{\omega} \leq \frac{1}{4}.
\]

The last step is to to verify the other condition that
$\bar{B}(x^0,\rho_{-}) \subseteq \bar{B}(x^0,2\gamma(\theta)r)$.
Noting that $\sqrt{1-2h} \geq 1-2h$ for $0 \leq h \leq 1/2$, it is
seen that
\begin{eqnarray}
\rho_-(\eta,\hat{\omega}) & = & \frac{1-\sqrt{1-2\eta\hat{\omega}}}{\hat{\omega}} \nonumber \\
& \leq & 2 \eta \nonumber \\
& \leq & 2 (\omega r + 1 ) r \nonumber \\
& \leq & 2\gamma(\theta)r. \qedhere \nonumber
\end{eqnarray}
\end{proof}

Next results are concerned with the size of the patch satisfying the exclusion test.

% Lemma 1.1: If ||f(x0)|| is small and roots are well conditioned and such, then x0 is near a root.
\begin{lemma}
\label{lem1}
Let $f: \mathbb{C}^n\rightarrow \mathbb{C}^n$ be a polynomial function with generic coefficients.
Assume that all zeros of $f$ are invertible.
Let $x^0$ be a point in $\mathbb{R}^n$.
If
\begin{equation}
\label{lemmaas}
\norm{\inv{f'(x^*)}f(x^0)} \leq \frac{1}{2\omega_f}
\end{equation}
for all complex zeros $x^*$ of $f$,
then there exists $\hat{x}^*$, a zero of $f$, such that
\begin{eqnarray}
\label{lemmacn}
\norm{x^0-\hat{x}^*} & \leq & \frac{1-\sqrt{1-2\omega_f\norm{\inv{f'(\hat{x}^*)}f(x^0)}}}{\omega_f} \\
& \equiv & \sigma(\hat{x}^*,x^0). \nonumber
\end{eqnarray}
\end{lemma}
\begin{proof}
By the assumption that $f$ has generic coefficients, the polynomial $f$ has a finite number of zeros.
Let $x^*_1$, $x^*_2$, \ldots, $x^*_d$ be all the complex zeros of $f$.
Recall that a multiple zero has singular Jacobian. Hence, $f$ has no multiple zeros by assumption.

Define the polynomial $\bar{f}(x) = f(x)-f(x^0)$.  Note that $x^0$ is a zero of $\bar{f}$.
We apply the Kantorovich's theorem for complex functions (see \cite{farouki}) to each $x^*_i$ with respect to
$\bar{f}$.  For each $x^*_i$, we use $D = \bar{B}(x^*_i,\rho_*(x^*_i))$ and $\omega = \omega_f$.
Since $\eta \equiv \norm{\inv{\bar{f}'(x^*_i)}\bar{f}(x^*_i)} = \norm{\inv{f'(x^*_i)}f(x^0)}$,
the assumption (\ref{lemmaas}) guarantees that the condition $\eta\omega \leq 1/2$ is satisfied.
The condition $\bar{B}(x^*_i,\rho_-) \subseteq D$ is also satisfied by the definition of $D$.  Therefore, the Kantorovich
theorem states that there is a zero of $\bar{f}$, call it $\bar{x}^*_i$, such that
\begin{equation}
\label{lenofh}
\norm{\bar{x}^*_i-x^*_i} \leq \sigma(x^*_i,x^0).
\end{equation}
Recall that, for any $j$, $x^*_j$ is the unique zero of $f$ in $\bar{B}(x^*_j,\rho_*(x^*_j))$.
Therefore,
\begin{equation}
\label{rootdistance}
\norm{x^*_i-x^*_j} > \max\{\rho_*(x^*_i),\rho_*(x^*_j)\}, i \neq j.
\end{equation}
But (\ref{lenofh}) and (\ref{rootdistance}) together imply that
\begin{equation}
\label{distinct}
\bar{x}^*_i \neq \bar{x}^*_j, i \neq j.
\end{equation}
Hence the mapping $x^*_i \rightarrow \bar{x}^*_i$ is injective.
But since $f$ has generic coefficients and $f$ and $\bar{f}$ are of the same degrees, $f$
has at least as many zeros as $\bar{f}$ \cite{emiris}.
This implies that $x^0 = \bar{x}^*_i$, for some $i$.  The lemma follows.
\end{proof}

\begin{thm}
\label{thm2} Let $f(x)=f(u,v)$ be a polynomial system in basis
$\phi_i(u)\phi_j(v)$ in two dimensions
with generic coefficients. %Let $m$ be the highest $u$ degree of
%$f$ and $n$ be its highest $v$ degree.
Let $x^0 = (u^0,v^0)$ be a
point in $[0,1]^2$ such that $f'(x^0)$ is invertible and $f(x^0)
\neq 0$, $x^*$ be the closest zero in $\mathbb{R}^2$ of $f$ to
$x^0$, and $\delta$ denote $\norm{x^0-x^*}$. Let $r > 0$ be such
that $\bar{B}(x^0,r) \subseteq [0,1]^2$. Assume $\delta >
\frac{1}{\omega_f}$.  Define $\hat{f}(\hat{u},\hat{v})$ such that
\begin{equation}
\label{rescale}
\hat{f}(\hat{u},\hat{v}) = f( \frac{2r}{h-l}\hat{u}-\frac{2hr}{h-l}+u^0+r,
\frac{2r}{h-l}\hat{v}-\frac{2hr}{h-l}+v^0+r).
\end{equation}
In other word, $\hat{f}$ is a polynomial in basis
$\phi_i(u)\phi_j(v)$ that reparametrizes with $[l,h]^2$ the
surface defined by $f$ over the patch $\bar{B}(x^0,r)$.
%Let $h(f)$ denote the convex hull of $\{\textbf{b}_{ij}\}$.
The bounding polygon of $\{ \hat{f}(u,v) : l \leq u,v \leq h\}$
satisfying item \ref{bounding_prop} of the basis properties listed in Section \ref{section_formulation}
does not contain the origin if
\begin{equation}
\label{deltahatcond} r \leq \frac{1}{2\theta\cond{f}^2}.
\end{equation}
\end{thm}
\begin{proof}
Let $X$ denote the patch $\bar{B}(x^0,r)$ and $x$ denote an arbitrary point in $X$.
Since $\delta > \frac{1}{\omega_f}$, the contrapositive of Lemma \ref{lem1} implies there exists
a zero $\bar{x}^*$ of $f$ satisfying
$\norm{f'(\bar{x}^*)^{-1}f(x^0)} > \frac{1}{2\omega_f}$. Therefore, the condition
(\ref{deltahatcond}) implies
\begin{eqnarray}
r & \leq & \frac{1}{2\theta\cond{f}^2} \nonumber \\
& \leq & \frac{1}{2\theta\omega_f\norm{f'(\bar{x}^*)^{-1}f'(x)}} \nonumber \\
& < &
\frac{\norm{f'(\bar{x}^*)^{-1}f(x^0)}}{\theta\norm{f'(\bar{x}^*)^{-1}f'(x)}}
\nonumber.
\end{eqnarray}
More specifically, we have
\[
r < \frac{\norm{f'(\bar{x}^*)^{-1}f(x^0)}}{\theta\max_{y\in
X}\norm{f'(\bar{x}^*)^{-1}f'(y)}},
\]
which is equivalent to
\begin{equation}
\label{fx1} \theta\cdot\max_{y \in
X}\norm{f'(\bar{x}^*)^{-1}f'(y)}\cdot r <
\norm{f'(\bar{x}^*)^{-1}f(x^0)}.
\end{equation}
Recall that $\max_{y \in X}\norm{f'(\bar{x}^*)^{-1}f'(y)}$ is the Lipschitz constant for
$\inv{f'(\bar{x}^*)}f$ on $X$. Hence, for any $x \in X$,
\begin{eqnarray}
\norm{f'(\bar{x}^*)^{-1}f(x)-f'(\bar{x}^*)^{-1}f(x^0)} &\leq &
\max_{y \in X}\norm{f'(\bar{x}^*)^{-1}f'(y)} \cdot \norm{x-x^0} \nonumber \\
 & \leq & \max_{y \in X}\norm{f'(\bar{x}^*)^{-1}f'(y)} \cdot r.
\label{lip}
\end{eqnarray}
Combining (\ref{fx1}) and (\ref{lip}) gives
\[
\theta\cdot\norm{f'(\bar{x}^*)^{-1}f(x)-f'(\bar{x}^*)^{-1}f(x^0)} <
\norm{f'(\bar{x}^*)^{-1}f(x^0)},
\]
which is equivalent to
\[
\theta\cdot\norm{f'(\bar{x}^*)^{-1}\hat{f}(\hat{x})-f'(\bar{x}^*)^{-1}\hat{f}(\hat{x}^0)}
< \norm{f'(\bar{x}^*)^{-1}\hat{f}(\hat{x}^0)}
\]
for some $\hat{x} \in [l,h]^2$, where $\hat{x}$ is the rescaled $x$ and $\hat{x}^0$ is the rescaled $x^0$
according to (\ref{rescale}). In particular,
\begin{equation}
\label{beforebi} \theta\cdot\max_{\hat{x} \in
[l,h]^2}\norm{f'(\bar{x}^*)^{-1}\hat{f}(\hat{x})-f'(\bar{x}^*)^{-1}\hat{f}(\hat{x}^0)}
< \norm{f'(\bar{x}^*)^{-1}\hat{f}(\hat{x}^0)}.
\end{equation}
Let $h(\hat{x}) \equiv f'(\bar{x}^*)^{-1}\hat{f}(\hat{x})$ and $g(\hat{x}) \equiv h(\hat{x})-h(\hat{x}^0)$. By (\ref{thetadef}),
\begin{equation}
\label{fromlemma} \norm{z} \leq
\theta\cdot\max_{\hat{x} \in
[l,h]^2}\norm{g(\hat{x})},
%\label{fromlemma} \norm{y-h(\hat{x}^0)} \leq
%2\tau\cdot\max_{\hat{x} \in
%[0,1]^2}\norm{h(\hat{x})-h(\hat{x}^0))},
\end{equation}
for any $z$ in the bounding polygon $P_g$ of $\{ g(\hat{x}) :
\hat{x} \in [l,h]^2 \}$. Since the bounding polygon is required to be
translationally invariant (item \ref{bounding_prop} of the basis properties
listed in Section \ref{section_formulation}),
(\ref{fromlemma}) is equivalent to
\begin{equation}
\label{shifted_theta} \norm{y-h(\hat{x}^0)} \leq
\theta\cdot\max_{\hat{x} \in
[l,h]^2}\norm{h(\hat{x})-h(\hat{x}^0))},
\end{equation}
for any $y$ in the bounding polygon $P_h$ of $\{ h(\hat{x}) :
\hat{x} \in [l,h]^2 \}$.
Substituting (\ref{shifted_theta}) into the
left hand side of (\ref{beforebi}) yields
\begin{equation}
\label{convexnozero} \norm{y-h(\hat{x}^0)} < \norm{h(\hat{x}^0)},
\end{equation}
which implies that $P$ does not contain the origin. Since
$f'(\bar{x}^*)^{-1}$ is invertible and the bounding polygon is
affinely invariant, the bounding polygon of $\{ \hat{f}(\hat{x})
: \hat{x} \in [l,h]^2 \}$ does not contain the origin, either.
\end{proof}

\begin{thm}
\label{thm3} Let $f(x)=f(u,v)$ be a polynomial system in basis
$\phi_i(u)\phi_j(v)$ in two dimensions with generic coefficients
whose zeros are sought.  Let $X = \bar{B}(x^0,r)$ be a patch under
consideration during the course of the KTS algorithm. The
algorithm does not need to subdivide $X$ if
\begin{equation}
\label{thm3ass} r \leq \frac{1}{2}\cdot\min
\left\{\frac{1-1/\gamma(\theta)}{\omega_f},
\frac{1}{2\theta\cond{f}^2}\right\}.
\end{equation}
\end{thm}
\begin{proof}
%If $X$ contains a zero $x^*$, (\ref{thm3ass}) implies that $\bar{B}(x^0,r)$ passes the Kantorovich test
%according to Theorem \ref{thm1}.
If $\delta > 1/\omega_f$, where $\delta$ is the distance between
$x^0$ and the closest zero $x^*$, $r \leq
\left(1-1/\gamma(\theta)\right)/(2\omega_f)$ implies that $X$ does
not contain a zero.  Therefore, $r \leq 1/(4\theta\cond{f}^2)$
implies that $X$ is excluded by the exclusion test according to
Theorem \ref{thm2}.

Observe that $\omega_* \leq \omega_f$. If $\delta \leq
\frac{1}{\omega_f}$, for any $x \in X$,
\begin{eqnarray}
\norm{x-x^*} & \leq & \norm{x-x^0} + \norm{x^0-x^*} \nonumber \\
& \leq & r + \delta \nonumber \\
& \leq & \frac{1-1/\gamma(\theta)}{\omega_f} + \frac{1}{\omega_f} \nonumber \\
& < & \frac{2}{\omega_f}  \leq \frac{2}{\omega_*} = \rho_*. \nonumber
\end{eqnarray}
In other word, $X$ is contained within $\bar{B}(x^*,\rho_*)$, a
safe region and therefore is excluded, provided that $x^*$ is
found before $X$ is checked against all safe regions. By Theorem
\ref{thm1}, $x^*$ is found by a region of size $2r \leq
(1-1/\gamma(\theta))/\omega_{D'}$. Since KTS examines larger
regions before smaller ones, $x^*$ is found before $X$ is checked
against safe regions.
\end{proof}
It should be noted that
\begin{equation}
\label{lingam}
1-1/\gamma(\theta) \geq 1/(18\theta)
\end{equation}
hence both terms of the
right hand side of (\ref{thm3ass}) are asymptotically linear in $1/\theta$.
The inequality (\ref{lingam}) follows from the fact that
\begin{equation}
\label{asmp}
\sqrt{1+a} \leq 1+a/2-a^2/9
\end{equation}
for any $a \in [0,1/4]$. To prove (\ref{asmp}), simplify (\ref{asmp}) to $a^2-9a+9/4 \geq 0$,
whose left hand side is a convex quadratic polynomial that crosses the x-axis at $(9-6\sqrt{2})/2 \approx .2574$
and at $(9+6\sqrt{2})/2$.

\section{Computational results}

The KTS algorithm is implemented in Matlab and is tested against a
number of problem instances with varying condition numbers.  As
B\'{e}zier surfaces are widely used in geometric modeling, we
choose to implement KTS for the Bernstein basis case.  Most of the
test problems are created by using normally distributed random
numbers as the coefficients $c_{ij}$'s of $f$.  For some of the
test problems especially those with high condition number, some
coefficients are manually entered.  The degrees of the test
polynomials are between biquadratic and biquartic.  As an example,
the test case with $\cond{f} = 3.5 \times 10^3$ is $c_{00} =
(1.2,.5)^T$, $c_{01} = (-.6,-.6)^T$, $c_{02} = (.1,1.1)^T$,
$c_{10} = (-1.1,-.3)^T$, $c_{11} = (.6,-2.3)^T$, $c_{12} =
(-2,-.1)^T$, $c_{20} = (.6,1.2)^T$, $c_{21} = (-1.1,-1.2)^T$, and
$c_{22} = (-.5,.4)^T$. This is the test problem for the result in
the second row of Table \ref{table_res}.

For the experiment, we use the algorithm by J\'{o}nsson and
Vavasis \cite{jonsson} to compute the complex zeros required to
estimate the condition number. Table \ref{table_res} compares the
efficiency of KTS with its condition number. The total number of
subpatches examined by KTS during the entire computation, the
width of the smallest patch among those examined, and the maximum
number of Newton iterations (in the cases with more than one zero)
to converge to a zero are reported.  The result shows that KTS
needs to examine more number of patches and needs to subdivide to
smaller patches as the condition number becomes larger.  Note that
the high number of Newton iterations of some test cases is due to
roundoff error.

%\begin{figure}
% \centering
% \includegraphics[width=0.9\textwidth]{fig2}
% \caption{Intersection between the chosen bicubic surface and the line of the fourth test case.}
% \label{fig}
%\end{figure}
\begin{table}
\begin{center}
\begin{tabular}{|r|r|r|r|r|r|}
\hline
& \multicolumn{1}{c}{Num. of}\vline & \multicolumn{1}{c}{Distance}\vline & \multicolumn{1}{c}{Num. of}\vline &
Smallest& \multicolumn{1}{c}{Max. num. of}\vline\\
$\cond{f}$& \multicolumn{1}{c}{zeros}\vline &  \multicolumn{1}{c}{between two}\vline &
 \multicolumn{1}{c}{patches}\vline &  \multicolumn{1}{c}{width}\vline &  \multicolumn{1}{c}{Newton}\vline \\
&&  \multicolumn{1}{c}{closest zeros}\vline  &examined& &iterations \\
\hline
$6.0 \times 10^2$ & 1 & - & 21 & .0625 & 3 \\
$3.5 \times 10^3$ & 2 & .4196 & 29 & .0625 & 3 \\
$8.3 \times 10^4$ & 2 & .6638 & 33 & .0625 & 3 \\
$1.6 \times 10^5$ & 1 & - & 41 & .03125 & 4 \\
$2.2 \times 10^7$ & 3 & .3624 & 57 & .03125 & 4 \\
$1.3 \times 10^8$ & 4 & .2806 & 81 & .015625 & 6 \\
$1.9 \times 10^9$ & 4 & .3069 & 69 & .03125 & 6 \\
$2.0 \times 10^{10}$ & 2 & .7810 & 105 & .015625 & 6 \\
$2.9 \times 10^{11}$ & 1 & - & 257 & .0039 & 9 \\
\hline
\end{tabular}
\end{center}
%\caption{The condition number, the number of zeros in $[0,1]^2$, and
%the distance between the closest pair of zeros in $[0,1]^2$ of each test problem instance. }
\caption{Efficiency of KTS algorithm on problems of different condition numbers.}
\label{table_res}
\end{table}

\section{Conclusion and future directions}

We present KTS algorithm for finding the intersections between a
parametric surface and a line. By using the combination of
subdivision and Kantorovich's theorem, our algorithm can take
advantage of the quadratic convergence of Newton's method without
the problems of divergence and missing some intersections that
commonly occur with Newton's method. We also show that the
efficiency of KTS has an upper bound that depends solely on the conditioning of the
problem and the representation of the surface. Nevertheless, there
are a number of questions left unanswered by this article such as
\begin{itemize}
\item \textbf{Extensibility to piecewise polynomial surfaces
and/or NURBS}. Since KTS only requires the ability to compute the
bounding polygon of a subpatch that satisfies the list of basis properties, it may
be possible to extend KTS to handle these more general surfaces if
bounding polygons having similar properties can be computed
relatively quickly.

\item \textbf{Tighter condition number}. The
condition number presented earlier seems overly loose.  It is
likely that a tighter condition number exists. If a tighter
condition number is found, we would be able to calculate a tighter
bound on the time complexity of KTS, too.

\item \textbf{The
necessity of the generic coefficients assumption}. Is it possible
to analyze the efficiency of KTS without this assumption?

\item
\textbf{Using KTS in floating point arithmetic}. In the presence
of roundoff error, we may need to make adjustments for KTS to be
able to guarantee that all zeros are found.  In addition, the
accuracy of the computed zeros would become an important issue in
this case.

\item \textbf{Choice of polynomial basis}. It is evident from
Table \ref{table_theta} that the Chebyshev basis has the best
(smallest) value of $\theta$, and therefore ought to require the
fewest number of subdivisions.  Our preliminary computational
results \cite{srijuntongsiri2} comparing bases, however, do not
indicate a clear-cut advantage for the Chebyshev basis. Therefore,
the impact of the choice of basis on practical efficiency is an
interesting topic for further research.
\end{itemize}

\end{document}